\documentclass[a4paper,10pt]{article}

\usepackage[utf8x]{inputenc}
\usepackage{amsmath,amssymb}
\newtheorem{thm}{Theorem}[section]
\newtheorem{lem}[thm]{Lemma}
\newtheorem{pro}[thm]{Proposition}
\newtheorem{cor}[thm]{Corollary}
\usepackage[all]{xy}

\title{Approximate Hermitian-Yang-Mills structures and semistability for Higgs bundles. \\
       I: Generalities and the one-dimensional case.}
\author{S. A. H. Cardona\footnote{Electronic address: sholguin@sissa.it}\\SISSA - Via Bonomea 265 - 34136, Trieste - Italy}

\begin{document}

\maketitle

\begin{abstract}
We review the notions of (weak) Hermitian-Yang-Mills structure and approximate Hermitian-Yang-Mills structure for Higgs bundles. Then, 
we construct the Donaldson functional for Higgs bundles over compact K\"ahler manifolds and we present some basic properties of it. In particular, 
we show that its gradient flow can be written in terms of the mean curvature of the Hitchin-Simpson connection. We also study some properties of 
the solutions of the evolution equation associated with that functional. Next, we study the problem of the existence of approximate Hermitian-Yang-Mills 
structures and its relation with the algebro-geometric notion of semistability and we show that for a compact Riemann surface, the 
notion of approximate Hermitian-Yang-Mills structure is in fact the differential-geometric counterpart of the notion of semistability. Finally, 
we review the notion of admissible Hermitian structure on a torsion-free Higgs sheaf and define the Donaldson functional for such an object.      
\end{abstract}

\section{Introduction}

In complex geometry, the Hitchin-Kobayashi correspondence asserts that the notion of (Mumford-Takemoto) stability, originally introduced in 
algebraic geometry, has a differential-geometric equivalent in terms of special metrics. In its classical version, this correspondence is 
established for holomorphic vector bundles over compact K\"ahler manifolds and says that such bundles are polystable if and only if they 
admit an Hermitian-Einstein\footnote{In the literature Hermitian-Einstein, Einstein-Hermite and Hermitian-Yang-Mills are all synonymous. Sometimes, 
also the terminology Hermitian-Yang-Mills-Higgs is used \cite{Biswas-Schumacher}.} structure. This correspondence is also true 
for Higgs bundles. \\
      
The history of this correspondence starts in 1965, when Narasimhan and Seshadri \cite{Narasimhan-Seshadri} proved that a holomorphic bundle 
on a Riemann surface is stable if and only if it corresponds to a projective irreducible representation of the fundamental group of the surface. 
Then, in the 80's Kobayashi \cite{1art-Kobayashi} introduced for the first time the notion of Hermitian-Einstein structure in a 
holomorphic vector bundle, as a generalization of a K\"ahler-Einstein metric in the tangent bundle. Shortly after, Kobayashi \cite{2art-Kobayashi} 
and L\"ubke \cite{Lubke 2} proved that a bundle with an irreducible Hermitian-Einstein structure must be necessarily stable. Donaldson in \cite{Donaldson-1} 
showed that the result of Narasimhan and Seshadri \cite{Narasimhan-Seshadri} can be formulated in terms of metrics and showed that the concepts of 
polystability and Hermitian-Einstein metrics are equivalent for holomorphic vector bundles over a compact Riemann surface. Then, Kobayashi and Hitchin conjectured 
that the equivalence should be true in general for holomorphic vector bundles over K\"ahler manifolds. However, the route starting from stability and 
showing the existence of special structures in higher dimensions took some time. \\

The existence of Hermitian-Einstein structures in a stable holomorphic vector bundle was proved by Donaldson for projective algebraic surfaces in 
\cite{Donaldson-2} and for projective algebraic manifolds in \cite{Donaldson-3}. Finally, Uhlenbeck and Yau showed this for general compact K\"ahler 
manifolds in \cite{Uhlenbeck-Yau} using some techniques from analysis and Yang-Mills theory. Hitchin \cite{Hitchin}, while studying the self-duality 
equations over a compact Riemann surface, introduced the notion of Higgs field, and showed that the result of Donaldson for Riemann surfaces could be 
modified to include the presence of a Higgs field. Following the results of Hitchin, Simpson in \cite{Simpson} defined a Higgs bundle to be a holomorphic 
bundle together with a Higgs field and proved the Hitchin-Kobayashi correspondence for such an object. As an application of this correspondence, 
Simpson \cite{Simpson}, \cite{Simpson 2} studied in detail a one-to-one correspondence between stable Higgs bundles over compact K\"ahler manifolds 
with vanishing Chern classes and representations of the fundamental group of the K\"ahler manifold. \\

The Hitchin-Kobayashi correspondence has been further extended in several directions. Simpson \cite{Simpson} studied the Higgs case for non-compact K\"ahler 
manifolds satisfying some additional requirements and L\"ubke and Teleman \cite{Lubke} studied the correspondence for compact complex manifolds. Bando 
and Siu \cite{Bando-Siu} extended this correspondence for torsion-free sheaves over compact K\"ahler manifolds and introduced the notion of admissible 
Hermitian metric for such objects. Following the ideas of Bando and Siu, Biswas and Schumacher \cite{Biswas-Schumacher} generalized this to the Higgs case. \\ 
  
In \cite{Donaldson-2} and \cite{Donaldson-3}, Donaldson introduced for the first time a functional (which is known as the Donaldson functional) and later 
Simpson \cite{Simpson} defined this functional for Higgs bundles. Kobayashi in \cite{Kobayashi} constructed the same functional 
in a different form and showed that it played a fundamental role in a possible extension of the Hitchin-Kobayashi correspondence. In fact, using that functional 
he proved that for holomorphic vector bundles over projective algebraic manifolds, the counterpart of semistability is the notion 
of approximate Hermitian-Einstein structure.  \\        
     
In this article we show that for a Higgs bundle on a compact Riemann surface, there is a correspondence between semistability and the existence of approximate 
Hermitian-Yang-Mills structures. We do this by using a Donaldson functional approach analogous to that of Kobayashi \cite{Kobayashi}. This 
result covers the classical case if we take the Higgs field equal to zero. \\ 

The correspondence between semistability and the existence of approximate Hermitian-Einstein metrics in the ordinary case has been studied 
recently in \cite{Adam} using a technique developed by Buchdahl \cite{Buchdahl-0} for the desigularization of sheaves in the case of compact complex surfaces. 
In a future work, following \cite{Biswas-Schumacher}, we will study in detail the notion of admissible Hermitian metric on a Higgs sheaf, the Donaldson 
functional for torsion-free sheaves and the correspondence between semistability and the existence of approximate Hermitian-Yang-Mills metrics for 
Higgs bundles on compact K\"ahler manifolds of higher dimensions. \\

\noindent{\bf Acknowledgements} \\

\noindent First to all, the author would like to thank his thesis advisor, Prof. U. Bruzzo, for his constant support and encouragement and also for suggesting 
this problem. I would like to thank also The International School for Advanced Studies (SISSA) for graduate fellowship support. Finally, I would like 
to thank O. Iena, G. Dossena, P. Giavedoni and Carlos Marin for many enlightening discussions and for some helpful comments. The 
results presented in this article will be part of my Ph.D. thesis in the Mathematical-Physics sector at SISSA.

\section{Preliminaries}
We start with some basic definitions. Let $X$ be an $n$-dimensional compact K\"ahler manifold with $\omega$ its K\"ahler form and let 
$\Omega_{X}^{1}$ be the cotangent sheaf to $X$, i.e., it is the sheaf of holomorphic one-forms on $X$. A Higgs sheaf ${\mathfrak E}$ over $X$ is a coherent sheaf $E$ over $X$, together with a 
morphism $\phi : E\rightarrow E\otimes\Omega_{X}^{1}$ of ${\cal O}_{X}$-modules, such that the morphism  
$\phi\wedge\phi : E\rightarrow E\otimes\Omega_{X}^{2}$ vanishes. The morphism $\phi$ is called the Higgs field of $\mathfrak E$. A Higgs 
sheaf $\mathfrak E$ is said to be torsion-free (resp. reflexive, locally free, normal, torsion) if the sheaf $E$ is torsion-free 
(resp. reflexive, locally free, normal, torsion). A Higgs subsheaf ${\mathfrak F}$ of ${\mathfrak E}$ is a subsheaf $F$ of $E$ such that $\phi(F)\subset F\otimes\Omega_{X}^{1}$. 
A Higgs bundle ${\mathfrak E}$ is just a Higgs sheaf in which the sheaf $E$ is a locally free ${\cal O}_{X}$-module. \\

Let ${\mathfrak E}_{1}$ and ${\mathfrak E}_{2}$ be two Higgs sheaves over a compact K\"ahler manifold $X$. A morphism 
between ${\mathfrak E}_{1}$ and ${\mathfrak E}_{2}$ is a map $f:E_{1}\longrightarrow E_{2}$ such that the diagram 
\begin{displaymath}
 \xymatrix{
E_{1} \ar[r]^{\phi_{1}} \ar[d]^{f}    &     E_{1}\otimes\Omega_{X}^{1} \ar[d]^{f\otimes{\rm Id}}   \\
E_{2} \ar[r]^{\phi_{2}}               &     E_{2}\otimes\Omega_{X}^{1}   \\
}
\end{displaymath}
commutes. We will denote such a morphism by $f:{\mathfrak E}_{1}\longrightarrow{\mathfrak E}_{2}$. A sequence of Higgs sheaves is a sequence of their 
corresponding coherent sheaves where each map is a morphism of Higgs sheaves. A short exact sequence of Higgs sheaves (also called 
an extension of Higgs sheaves or a Higgs extension \cite{Simpson 2}, \cite{Bradlow-Gomez}) is defined in the obvious way. \\    

We define the degree ${\rm deg\,}{\mathfrak E}$ and rank ${\rm rk\,}{\mathfrak E}$ of a Higgs sheaf simply as the degree and rank of the sheaf $E$. If the 
rank is positive, we introduce the quotient $\mu({\mathfrak E})={\rm deg\,}{\mathfrak E}/{\rm rk\,}{\mathfrak E}$ and call it the slope of the Higgs sheaf 
${\mathfrak E}$. In a similar way as in the ordinary case (see for instance \cite{Kobayashi}, \cite{Uhlenbeck-Yau}, \cite{Siu}, \cite{Bando-Siu}) there is a notion of stability 
for Higgs sheaves, which depends on the K\"ahler form and makes reference only to Higgs subsheaves \cite{Simpson}, \cite{Simpson 2}, \cite{Bruzzo-Granha}, 
\cite{Biswas-Schumacher}, \cite{Bradlow-Gomez}. 
Namely, a Higgs sheaf ${\mathfrak E}$ is said to be $\omega$-stable (resp. $\omega$-semistable) if it is torsion-free and for any Higgs subsheaf 
${\mathfrak F}$ with $0<{\rm rk\,}{\mathfrak F}<{\rm rk\,}{\mathfrak E}$ and torsion-free quotient, one has the inequality 
$\mu({\mathfrak F})<\mu(\mathfrak E)$ (resp. $\le$). We say that a Higgs sheaf is $\omega$-polystable if it decomposes into a direct sum 
of $\omega$-stable Higgs sheaves all having the same slope. \\

Let $\mathfrak E =(E,\phi)$ be a Higgs bundle of rank $r$ over $X$ and let $\omega$ be the K\"ahler form of $X$. Using the Chern connection 
$D_{(E,h)}$ of $E$ and the Higgs field $\phi$ one defines the Hitchin-Simpson connection on $\mathfrak E$ by:
\begin{equation}
{\mathcal D}_{({\mathfrak E},h)} = D_{(E,h)} + \phi + \bar\phi_{h}\,,
\end{equation}
where $\bar\phi_{h}$ is the adjoint of the Higgs field with respect to the Hermitian structure $h$, that is, it is defined 
by the formula $h(\bar\phi_{h}s,s')=h(s,\phi s')$ with $s,s'$ sections of the Higgs bundle. The curvature of the Hitchin-Simpson connection is then 
given by ${\mathcal R}_{(\mathfrak E,h)}={\mathcal D}_{(\mathfrak E,h)} \circ {\mathcal D}_{(\mathfrak E,h)}$ and hence
\begin{equation}
{\mathcal R}_{(\mathfrak E,h)} = R_{(E,h)} +  D'_{(E,h)}(\phi) +  D''_{(E,h)}(\bar\phi_{h}) + [\phi,\bar\phi_{h}]\,.   \label{Hitchin-Simpson curv}
\end{equation}
We say that the pair $({\mathfrak E},h)$ is Hermitian flat, if the curvature ${\mathcal R}_{(\mathfrak E,h)}$ vanishes. We denote by Herm$({\mathfrak E})$ the 
space of Hermitian forms in ${\mathfrak E}$ and by Herm$^{+}({\mathfrak E})$ the space of Hermitian structures (i.e., positive definite Hermitian forms) 
in ${\mathfrak E}$. For any Hermitian structure $h$ it is possible to identify ${\rm Herm}(\mathfrak E)$ with the tangent space of 
${\rm Herm}^{+}(\mathfrak E)$ at that $h$ (see \cite{Kobayashi} for details). That is
\begin{equation}
{\rm Herm}(\mathfrak E) = T_{h}\, {\rm Herm}^{+}(\mathfrak E)\,.
\end{equation}
If $v$ denotes an element in Herm$({\mathfrak E})$, one defines the endomorphism $h^{-1}v$ of 
${\mathfrak E}$ by the formula 
\begin{equation}
v(s,s')=h(s,h^{-1}vs'), \label{K endomorphism vs. form}
\end{equation}
where $s,s'$ are sections of ${\mathfrak E}$. We define a Riemannian structure in ${\rm Herm}^{+}(\mathfrak E)$ via this identification. Namely, for any $v,v'$ in 
${\rm Herm}(\mathfrak E)$ we define
\begin{equation}
(v,v')_{h}= \int_{X}{\rm tr}(h^{-1}v\cdot h^{-1}v')\,\omega^{n}/n!\,.   \label{inner product}
\end{equation}

The Higgs field $\phi$ can be considered as a section of ${\rm End}E\otimes\Omega^{1}_{X}$ and hence we have a natural dual morphism 
$\phi^{*}:E^{*}\rightarrow E^{*}\otimes\Omega^{1}_{X}$. From this it follows that ${\mathfrak E}^{*}=(E^{*},\phi^{*})$ is a Higgs bundle. On the other 
hand, if $Y$ is another compact K\"ahler manifold and $f:Y\rightarrow X$ is a holomorphic map, the pair defined by 
$f^{*}\mathfrak E= (f^{*}E,f^{*}\phi)$ is also a Higgs bundle. We have also some natural properties associated to tensor products and direct sums. 
In particular we have

\begin{pro}\label{basic prop. HYM}
Let ${\mathfrak E}_{1}$ and  ${\mathfrak E}_{2}$ two Higgs bundles with Higgs fields $\phi_{1}$ and $\phi_{2}$ respectively. Then \\
\noindent{\bf (i)} The pair ${\mathfrak E}_{1}\otimes{\mathfrak E}_{2}=(E_{1}\otimes E_{2},\phi)$ is a Higgs bundle with 
$\phi=\phi_{1}\otimes I_{2} + I_{1}\otimes\phi_{2}$. \\
{\bf (ii)} If ${\rm pr}_{i}:E_{1}\oplus E_{2}\rightarrow E_{i}$ with $i=1,2$ denote the natural projections, then 
${\mathfrak E}_{1}\oplus {\mathfrak E}_{2}=(E_{1}\oplus E_{2},\phi)$ is a Higgs bundle with 
$\phi={\rm pr}^{*}_{1}\phi_{1} + {\rm pr}^{*}_{2}\phi_{2}$.     
\end{pro}

In a similar form as in the ordinary case \cite{Kobayashi}, \cite{Siu}, \cite{Uhlenbeck-Yau}, we have a notion of Hermitian-Einstein 
structure for Higgs bundles \cite{Simpson}, \cite{Simpson 2}. Let us consider the usual star operator $\ast : A^{p,q}\rightarrow A^{n-q,n-p}$ and 
the operator $L: A^{p,q}\rightarrow A^{p+1,q+1}$ defined by $L\varphi=\omega\wedge\varphi$, where $\varphi$ is a form on $X$ of type $(p,q)$. Then 
we define, as usual, 
$\Lambda=\ast^{-1}\circ L\circ\ast: A^{p,q}\rightarrow A^{p-1,q-1}$. Consider now a metric $h$ in ${\rm Herm}^{+}(\mathfrak E)$ (i.e., $h$ is an 
Hermitian structure on ${\mathfrak E}$), associated with this metric we have a Hitchin-Simpson curvature ${\cal R}_{(\mathfrak E,h)}$. We can 
define the mean curvature of the Hitchin-Simpson connection, just by contraction of this curvature with the operator $i\Lambda$. In other words, 
\begin{equation}
{\cal K}_{(\mathfrak E,h)}=i\Lambda{\cal R}_{(\mathfrak E,h)}\,.  
\end{equation}
The mean curvature is an endomorphism\footnote{If we consider a local frame field $\{ e_{i}\}_{i=1}^{r}$ for $\mathfrak E$ and a local coordinate 
system $\{ z_{\alpha}\}_{\alpha=1}^{n}$ of $X$, the components of the mean curvature are given by ${\cal K}^{i}_{\,j}=
\omega^{\alpha\bar\beta}{\cal R}^{i}_{\,j\alpha\bar\beta}$.} in ${\rm End}(E)$. We say that $h$ is a {\it weak Hermitian-Yang-Mills structure} with 
factor $\gamma$ for $\mathfrak E$ if 
\begin{equation}
{\cal K}_{(\mathfrak E,h)}=\gamma\cdot I 
\end{equation}
where $\gamma$ is a real function defined on $X$ and $I$ is the identity endomorphism on $E$. From this definition it follows that if $h$ is a weak Hermitian-Yang-Mills structure with factor $\gamma$ 
for $\mathfrak E$, then the dual metric $h^{\ast}$ is a weak Hermitian-Yang-Mills structure for the dual bundle ${\mathfrak E}^{\ast}$ and also, that any 
metric $h$ on a Higgs line bundle is necessarily a weak Hermitian-Yang-Mills structure. As in the ordinary case, also for Higgs bundles we have some 
simple properties related with the notion of weak Hermitian-Yang-Mills structure, in particular, from the usual formulas for the curvature of tensor 
products and direct sums we have the following

\begin{pro}\label{basic prop. HYM} {\bf (i)} If $h_{1}$ and $h_{2}$ are two weak Hermitian-Yang-Mills structures with factors $\gamma_{1}$ and  $\gamma_{2}$ for 
Higgs bundles ${\mathfrak E}_{1}$ and ${\mathfrak E}_{2}$, respectively, then $h_{1}\otimes h_{2}$ is a weak Hermitian-Yang-Mills structure with 
factor $\gamma_{1}+\gamma_{2}$ for the tensor product bundle ${\mathfrak E}_{1}\otimes {\mathfrak E}_{2}$. \\
{\bf (ii)} The metric $h_{1}\oplus h_{2}$ is a weak Hermitian-Yang-Mills structure with factor $\gamma$ for the Whitney sum 
${\mathfrak E}_{1}\oplus {\mathfrak E}_{2}$ if and only if both metrics $h_{1}$ and $h_{2}$ are weak Hermitian-Yang-Mills structures with the same 
factor $\gamma$ for ${\mathfrak E}_{1}$ and ${\mathfrak E}_{2}$, respectively.
\end{pro}

If we have a weak Hermitian-Yang-Mills structure in which the factor $\gamma=c$ is constant, we say that $h$ is an {\it Hermitian-Yang-Mills 
structure} with factor $c$ for $\mathfrak E$. From Proposition \ref{basic prop. HYM} and this definition we get 

\begin{cor}\label{} Let $h\in {\rm Herm}^{+}(\mathfrak E)$ be a (weak) Hermitian-Yang-Mills structure with factor $\gamma$ for the Higgs 
bundle $\mathfrak E$. Then \\
{\bf (i)} The induced Hermitian metric on the tensor product ${\mathfrak E}^{\otimes p}\otimes {\mathfrak E}^{*\,\otimes q}$ is a (weak) 
Hermitian-Yang-Mills structure with factor $(p-q)\gamma$. \\
{\bf (ii)} The induced Hermitian metric on $\bigwedge^{p}{\mathfrak E}$ is a (weak) Hermitian-Yang-Mills structure with factor $p\gamma$ for every 
$p\le r={\rm rk}{\mathfrak E}$.  
\end{cor}

In general, if $h$ is a weak Hermitian-Yang-Mills structure with factor $\gamma$, the slope of $\mathfrak E$ can be written in terms of $\gamma$. 
To be precise, we obtain 

\begin{pro}\label{} If $h\in {\rm Herm}^{+}(\mathfrak E)$ is a weak Hermitian-Yang-Mills structure with factor $\gamma$, then
\begin{equation}
{\mu}({\mathfrak E})=\frac{1}{2n\pi}\int_{X}\gamma\,\omega^{n}.    \label{c and degree}
\end{equation}
\end{pro}

\noindent {\it Proof:} Let ${\cal R}$ be the Hitchin-Simpson curvature of $\mathfrak E$, then in general we have the identity
\begin{equation}
in{\cal R}\wedge\omega^{n-1}={\cal K}\,\omega^{n}\,.    \label{general id}
\end{equation}
Now, by hypothesis $h$ is a weak Hermitian-Yang-Mills structure with factor $\gamma$, then taking the trace of (\ref{general id}) and integrating 
over $X$ we obtain\footnote{We consider here the integral $\frac{i}{2\pi}\int_{X}{\rm tr}\,{\mathcal R}\wedge\omega^{n-1}.$ 
Notice that only the $(1,1)$ part of the curvature ${\mathcal R}$ makes a real contribution in such an integral and since 
${\rm tr\,}[\phi,\bar\phi]$ is identically zero, that integral must be the degree of the holomorphic bundle $E$ and hence it is 
equal to ${\rm deg\,}{\mathfrak E}$.} 
\begin{equation}
{\rm deg\,}{\mathfrak E}= \frac{r}{2n\pi}\int_{X}\gamma\,\omega^{n}\,,
\end{equation}
where $r$ is the rank of ${\mathfrak E}$.  \;\;Q.E.D.  \\

Consider now a real positive function $a=a(x)$ on $X$, then $h'=ah$ defines another Hermitian metric on $\mathfrak E$. Since 
$h'$ is a conformal change of $h$, we have in particular $\bar\phi_{h'}=\bar\phi_{h}$. Then, from (\ref{general id}) we obtain
\begin{equation}
{\cal K}'\,\omega^{n} = in(R' + [\phi,\bar\phi_{h'}])\wedge\omega^{n-1} = \left(K' + i\Lambda[\phi,\bar\phi_{h}]\right)\omega^{n}. \label{id K's}
\end{equation}
Now, defining $\Box_{0}=i\Lambda d''d'$  (see \cite{Kobayashi} for details) we have $K'=K + \Box_{0}({\rm log}\,a)$ and hence using 
the identity (\ref{id K's}) we get
\begin{equation}
{\cal K}'\,\omega^{n} = {\cal K}\,\omega^{n} + \Box_{0}({\rm log}\,a)\,\omega^{n}.
\end{equation}
From this we conclude the following 

\begin{lem}\label{conformal change} Let $h$ be a weak Hermitian-Yang-Mills structure with factor $\gamma$ for $\mathfrak E$ and let $a$ be 
a real positive definite function on $X$, then $h'=ah$ is a weak Hermitian-Yang-Mills structure with factor $\gamma'=\gamma + \Box_{0}({\rm log}\,a)$. 
\end{lem}

Making use of Lemma \ref{conformal change} we can define a constant $c$ which plays an important role in the definition of the Donaldson functional. Such a 
constant $c$ is an average of the factor $\gamma$ of a weak Hermitian-Yang-Mills structure. Namely 

\begin{pro}\label{} If $h\in {\rm Herm}^{+}(\mathfrak E)$ is a weak Hermitian-Yang-Mills structure with factor $\gamma$, then there 
exists a conformal change $h'=ah$ such that $h'$ is an Hermitian-Yang-Mills structure with constant factor $c$, given by
\begin{equation}
c\int_{X}\omega^{n}=\int_{X}\gamma\,\omega^{n}\,. \label{def. c}
\end{equation}
Such a conformal change is unique up to homothety. 
\end{pro}
\noindent {\it Proof:} Let $c$ be as in (\ref{def. c}), then 
\begin{equation}
 \int_{X}(c-\gamma)\omega^{n}=0\,.   \label{cond. c}
\end{equation}
It is sufficient to prove that there is a function $u$ satisfying the equation 
\begin{equation}
 \Box_{0}u=c-\gamma\,, \label{Hodge eq.}
\end{equation}
where, as we said before $\Box_{0}=i\Lambda d''d'$. Because if this holds, then by applying Lemma \ref{conformal change} with the function $a=e^{u}$ the result 
follows. \\

Now, from Hodge theory we know that the equation (\ref{Hodge eq.}) has a solution if and only if the function $c-\gamma$ is orthogonal 
to all $\Box_{0}$-harmonic functions. Since $X$ is compact, a function is $\Box_{0}$-harmonic if and only if it is constant. But (\ref{cond. c}) says 
that $c-\gamma$ is orthogonal to the constant functions and hence the equation (\ref{Hodge eq.}) has a solution $u$. Finally, the uniqueness follows from 
the fact that $\Box_{0}$-harmonic functions are constant.  \;\;Q.E.D.  \\

Since every weak Hermitian-Yang-Mills structure can be transformed into an Hermitian-Yang-Mills structure using a conformal change of the metric, without 
loss of generality we avoid using weak structures and work directly with Hermitian-Yang-Mills structures.

\section{Approximate Hermitian-Yang-Mills structures}
As we have seen in the preceding section, if we have an Hermitian-Yang-Mills structure with factor $c$, this constant can be evaluated directly 
from (\ref{c and degree}) and we have
\begin{equation}
 c=\frac{2\pi\,\mu(\mathfrak E)}{(n-1)!\,{\rm vol}X}\,.   \label{general def. c}
\end{equation}
On the other hand, regardless if we have an Hermitian-Yang-Mills structure or not on $\mathfrak E$, we can always define a constant $c$ just by 
(\ref{general def. c}). Introduced in such a way, $c$ depends only on $c_{1}(\mathfrak E)$ and the cohomology class of $\omega$ and not on the metric 
$h$. We define the length of the endomorphism ${\cal K}-cI$ by the formula
\begin{equation}
|{\cal K}-cI|^{2}={\rm tr}\left[({\cal K}-cI)\circ({\cal K}-cI)\right]\,. 
\end{equation}
We say that a Higgs bundle $\mathfrak E$ over a compact K\"ahler manifold $X$ admits an {\it approximate Hermitian-Yang-Mills structure} 
if for any $\epsilon>0$ there exists a metric $h$ (which depends on $\epsilon$) such that
\begin{equation}
{\max}|{\cal K}-cI|<\epsilon\,.    \label{approx. HYM ineq}
\end{equation}

From the above definition it follows that ${\mathfrak E}^{*}$ admits an approximate Hermitian-Yang-Mills structure if $\mathfrak E$ does. This notion 
satisfies some simple properties with respect to tensor products and direct sums.

\begin{pro}\label{} If ${\mathfrak E}_{1}$ and ${\mathfrak E}_{2}$ admit approximate Hermitian-Yang-Mills structures, so does their tensor 
product ${\mathfrak E}_{1}\otimes{\mathfrak E}_{2}$. Furthermore if $\mu({\mathfrak E}_{1})=\mu({\mathfrak E}_{2})$, so does their Whitney sum 
${\mathfrak E}_{1}\oplus{\mathfrak E}_{2}$\,.
\end{pro}

\noindent{\it Proof:} Assume that ${\mathfrak E}_{1}$ and ${\mathfrak E}_{2}$ admit approximate Hermitian-Yang-Mills structures with factors 
$c_{1}$ and $c_{2}$ respectively and let $\epsilon>0$. Then, there exist $h_{1}$ and $h_{2}$ such that 
\begin{equation}
\max_{X}|{\cal K}_{1}-c_{1}I_{1}|<\frac{\epsilon}{2\sqrt r_{2}}\,,   \quad\quad    \max_{X}|{\cal K}_{2}-c_{2}I_{2}|< \frac{\epsilon}{2\sqrt r_{1}}\,, \nonumber
\end{equation}
where $r_{1}, r_{2}$ and $I_{1}, I_{2}$ are the ranks and the identity endomorphisms of ${\mathfrak E}_{1}$ and ${\mathfrak E}_{2}$ respectively. 
Now, let $\cal K$ be the Hitchin-Simpson mean curvature of ${\mathfrak E}_{1}\otimes{\mathfrak E}_{2}$ associated with the metric 
$h=h_{1}\otimes h_{2}$. Then, by defining $c=c_{1}+c_{2}$ and $I=I_{1}\otimes I_{2}$ it follows
\begin{eqnarray*}
 |{\cal K}-cI| &=&  |{\cal K}_{1}\otimes I_{2} + I_{1}\otimes{\cal K}_{2} - (c_{1}+c_{2})I_{I}\otimes I_{2}| \\
              &\le& |({\cal K}_{1}-c_{1}I_{1})\otimes I_{2}| + |I_{1}\otimes({\cal K}_{2}-c_{2}I_{2})| \\
              &\le& \sqrt{r_{2}}\,|{\cal K}_{1}-c_{1}I_{1}| + \sqrt{r_{1}}\,|{\cal K}_{2}-c_{2}I_{2}| \\
               &<&  \epsilon  
\end{eqnarray*}
and hence the tensor product ${\mathfrak E}_{1}\otimes{\mathfrak E}_{2}$ admits an approximate Hermitian-Yang-Mills structure.\\

On the other hand, if $\mu({\mathfrak E}_{1})=\mu({\mathfrak E}_{2})$, necesarily the constants $c_{1}$ and $c_{2}$ coincide. Then, taking
this time $c=c_{1}=c_{2}$, $I=I_{1}\oplus I_{2}$ and ${\cal K}={\cal K}_{1}\oplus{\cal K}_{2}$, we have
\begin{eqnarray*}
 |{\cal K}-cI| &=&  |{\cal K}_{1}\oplus{\cal K}_{2} - c\,I_{1}\oplus I_{2}|\\
               &=&  \sqrt{{\rm tr}\,({\cal K}_{1}-c_{1}I_{1})^{2} + {\rm tr}\,({\cal K}_{2}-c_{2}I_{2})^{2}}\\
              &\le& |{\cal K}_{1}-c_{1}I_{1}| +  |{\cal K}_{2}-c_{2}I_{2}|\,.  
\end{eqnarray*}
From this inequality it follows that ${\mathfrak E}_{1}\oplus{\mathfrak E}_{2}$ admits an approximate Hermitian-Yang-Mills structure.  
  \;\;Q.E.D.


\begin{cor}\label{} If $\mathfrak E$ admits an approximate Hermitian-Yang-Mills structure, so do the tensor product bundle
${\mathfrak E}^{\otimes p}\otimes{\mathfrak E}^{*\otimes q}$ and the exterior product bundle $\bigwedge^{p}{\mathfrak E}$ whenever $p\le r$. 
\end{cor}

Finally, in a similar way as in the classical case, we have a version of the Bogomolov-L\"ubke inequality also for Higgs bundles admiting an 
approximate Hermitian-Yang-Mills structure (see \cite{Kobayashi}, \cite{Lenna} for details). To be precise, we obtain 

\begin{thm}\label{Lubke ineq} Let $\mathfrak E$ be a Higgs bundle over a compact K\"ahler manifold $X$ and suppose that $\mathfrak E$ admits an approximate 
Hermitian-Yang-Mills structure, then
\begin{equation}
\int_{X}\left[2r\,c_{2}(\mathfrak E) - (r-1)c_{1}(\mathfrak E)^{2}\right]\wedge\omega^{n-2} \ge 0\,.  \label{Lubke inequality}
\end{equation}
\end{thm}

\noindent{\it Proof:} Assume that $\mathfrak E$ admits an approximate Hermitian-Yang-Mills structure. Let $\epsilon>0$ and suppose $h_{\epsilon}$ 
is a metric on $\mathfrak E$ satisfying (\ref{approx. HYM ineq}). Then, we have closed $2k$-forms $c_{k}(\mathfrak E,h_{\epsilon})$ representing 
the $k$-th Chern classes. From \cite{Kobayashi}, Ch.IV, we obtain
\begin{equation}
(2r\,c_{2}(\mathfrak E,h_{\epsilon}) - (r-1)\,c_{1}(\mathfrak E,h_{\epsilon})^{2})\wedge\frac{\omega^{n-2}}{(n-2)!} = \left[r(|R_{\epsilon}|^{2} - 
|K_{\epsilon}|^{2}) + \sigma_{\epsilon}^{2} - |\rho_{\epsilon}|^{2}\right]\frac{\omega^{n}}{n!}\, \nonumber
\end{equation}
where the quantities on the right-hand side are associated to the metric $h_{\epsilon}$ and are given by $|K_{\epsilon}|^{2}={\rm tr}\,K_{\epsilon}^{2}$ 
and $\sigma_{\epsilon} = {\rm tr}\,K_{\epsilon}$ and 
\begin{equation}
|R_{\epsilon}|^{2}=\sum_{i,j,\alpha,\beta}|(R_{\epsilon})^{i}_{j\alpha\bar\beta}|^{2}\,, \quad\quad 
|\rho_{\epsilon}|^{2}=\sum_{i,\alpha,\beta}|(R_{\epsilon})^{i}_{i\alpha\bar\beta}|^{2}\,.     \nonumber
\end{equation}
Now, one has $r|R_{\epsilon}|^{2} \ge |\rho_{\epsilon}|^{2}$, and hence integrating over $X$ we obtain 
\begin{equation}
\int_{X}(2r\,c_{2}(\mathfrak E,h_{\epsilon}) - (r-1)\,c_{1}(\mathfrak E,h_{\epsilon})^{2})\wedge\frac{\omega^{n-2}}{(n-2)!} \ge \int_{X}\left[\sigma_{\epsilon}^{2} 
-r\,|K_{\epsilon}|^{2}\right]\frac{\omega^{n}}{n!}\,.    \label{general ineq}
\end{equation}
Since $h_{\epsilon}$ is an approximate Hermitian-Yang-Mills structure, we have 
\begin{equation}
\epsilon^{2}>|{\mathcal K}_{\epsilon}-cI|^{2}= |{\mathcal K}_{\epsilon}|^{2} - 2c\,\sigma_{\epsilon} + c^{2}r\,.   \label{estimative K epsilon 1}
\end{equation}
On the other hand, 
\begin{eqnarray*}
 |{\cal K}_{\epsilon}|^{2} &=& {\rm tr}\left[(K_{\epsilon} + i\Lambda[\phi,\bar\phi_{\epsilon}])\cdot(K_{\epsilon} + i\Lambda[\phi,\bar\phi_{\epsilon}])\right] \\
                           &=& |K_{\epsilon}|^{2} + 2\,i\Lambda\,{\rm tr}\left[K_{\epsilon}\cdot[\phi,\bar\phi_{\epsilon}]\right] + (i\Lambda)^{2}{\rm tr}\left[[\phi,\bar\phi_{\epsilon}]^{2}\right] \\
                           &=& |K_{\epsilon}|^{2} + 2\,i\Lambda\,{\rm tr}\left[{\cal K}_{\epsilon}\cdot[\phi,\bar\phi_{\epsilon}]\right]\,. 
\end{eqnarray*}
Now, ${\cal K}_{\epsilon}=cI + \epsilon\,A$ with $A$ a self-adjoint endomorphism of $E$ and hence we can estimate the term involving the trace in 
the last expression as
\begin{equation}
{\rm tr}\left[{\cal K}_{\epsilon}\cdot[\phi,\bar\phi_{\epsilon}]\right] = c\,{\rm tr}\,[\phi,\bar\phi_{\epsilon}] + \epsilon\,{\rm tr}\left[A\cdot[\phi,\bar\phi_{\epsilon}]\right] = \epsilon\,\eta 
\end{equation}
where the $(1,1)$-form $\eta = {\rm tr}\left[A\cdot[\phi,\bar\phi_{\epsilon}]\right]$. Consequently
\begin{equation}
|{\cal K}_{\epsilon}|^{2} = |K_{\epsilon}|^{2} + 2\epsilon\,(i\Lambda\eta)\,. \label{estimative K epsilon 2} 
\end{equation}
Finally, from (\ref{estimative K epsilon 1}) and (\ref{estimative K epsilon 2}) it follows
\begin{equation}
 \sigma_{\epsilon}^{2} - r\,|K_{\epsilon}|^{2} > (\sigma_{\epsilon} - cr)^{2} + f(\epsilon)  \nonumber
\end{equation}
where $f(\epsilon)=r\epsilon(2\,(i\Lambda\eta) - \epsilon)$. Then, by replacing this last expression in (\ref{general ineq}) we conclude
\begin{equation}
\int_{X}(2r\,c_{2}(\mathfrak E,h_{\epsilon}) - (r-1)\,c_{1}(\mathfrak E,h_{\epsilon})^{2})\wedge\frac{\omega^{n-2}}{(n-2)!} > \int_{X} f(\epsilon)\,\frac{\omega^{n}}{n!}\,.  \label{general ineq 2}
\end{equation}
Now, the integral on the left-hand side is independent of the metric $h_{\epsilon}$. On the other hand, the above inequality holds for all 
$\epsilon>0$ and clearly $f(\epsilon)\rightarrow 0$ as $\epsilon\rightarrow 0$. Therefore, one has the inequality 
(\ref{Lubke inequality}) if $\mathfrak E$ admits an approximate Hermitian-Yang-Mills metric.  \;\;Q.E.D. \\


We want to construct a functional $\mathcal L$ on Herm$^{+}({\mathfrak E})$, that will be called Donaldson's functional and whose gradient is related 
with the mean curvature of the Hitchin-Simpson connection. The construction of this functional is in certain way similar to the ordinary case. However, 
there are some differences, which in essence are due to the extra terms involving the Higgs field $\phi$ in the expression for the curvature 
(\ref{Hitchin-Simpson curv}).

\section{Donaldson's functional}

Given two Hermitian structures $h,k$ in Herm$^{+}({\mathfrak E})$, we connect them by a curve $h_{t}$, $0\le t \le 1$, in 
Herm$^{+}({\mathfrak E})$ so that $k=h_{0}$ and $h=h_{1}$. We set 
\begin{equation}
Q_{1}(h,k) = {\rm log}({\rm det}(k^{-1}h))\,, \quad\quad Q_{2}(h,k) = i \int_{0}^{1} {\rm tr} (v_{t}\cdot {\mathcal R}_{t})\,dt\,,
\end{equation}
where $v_{t}=h_{t}^{-1}\partial_{t}h_{t}$ and ${\mathcal R}_{t}$ denotes the curvature of the Hitchin-Simpson connection associated with $h_{t}$. Notice 
that $Q_{1}(h,k)$ does not involve the curve (in fact, it is the same functional of the ordinary case). On the other hand, the definition of 
$Q_{2}(h,k)$ uses explicitly the curve and differs from the ordinary case because of the extra terms in (\ref{Hitchin-Simpson curv}). We define
 the Donaldson functional by  
\begin{equation}
{\mathcal L}(h,k)= \int_{X}\left[Q_{2}(h,k) - \frac{c}{n} Q_{1}(h,k)\,\omega \right]\wedge \omega^{n-1}/(n-1)!\,,
\end{equation}
where $c$ is the constant given by 
\begin{equation}
c=\frac{2\pi\,\mu(\mathfrak E)}{(n-1)!\,{\rm vol\,}X}\,.
\end{equation}

Notice that the components of $(2,0)$ and $(0,2)$ type of ${\mathcal R}_{t}$ do not contribute to ${\mathcal L}(h,k)$.
This means that, in practice, it is enough to consider in the definition of $Q_{2}(h,k)$ just components of 
$(1,1)$-type\footnote{In other words, in computations involving integration over $X$, we can always replace the curvature 
by ${\mathcal R}_{t}^ {1,1}= R_{t} + [\phi,\bar\phi_{h_{t}}]$.}. \\

The following Lemma and the subsequent Proposition are straightforward generalizations of a result of Kobayashi 
(see \cite{Kobayashi}, Ch.VI, Lemma 3.6) to the Higgs case. Part of the proof is similar to the proof presented in \cite{Kobayashi}, 
however some differences arise because of the term involving the commutator in the Hitchin-Simpson curvature. 

\begin{lem}\label{Main lemma L} Let $h_{t}, a\le t \le b\,,$ be any differentiable curve in {\rm Herm}$^{+}({\mathfrak E})$ and $k$ any fixed Hermitian 
structure of ${\mathfrak E}$. Then, the (1,1)-component of
\begin{equation}
i \int_{a}^{b} {\rm tr} (v_{t}\cdot {\mathcal R}_{t})\,dt  +  Q_{2}(h_{a},k) - Q_{2}(h_{b},k)                     \label{first Lemma} 
\end{equation}  
is an element in $d'A^{0,1} + d''A^{1,0}$. 
\end{lem}

\noindent{\it Proof:} Following \cite{Kobayashi}, we consider the domain  $\Delta$ in ${\mathbb R}^{2}$ defined by  
\begin{equation}
\Delta = \{(t,s)| a\le t\le b\,, 0\le s \le 1\}
\end{equation}
and let $h:\Delta \rightarrow {\rm Herm}^{+}({\mathfrak E})$ be a smooth mapping such that $h(t,0)=k\,, h(t,1)=h_{t}$ for 
$a\le t \le b$\,, let $h(a,s)$ and $h(b,s)$ line segments curves\footnote{Notice that we have a simple expression for line segments curves from 
$k$ to $h_{t}$ given by $h(t,s)=sh_{t} + (1-s)k$.} from $k$ to $h_{a}$ and respectively from $k$ to $h_{b}$. We 
define the endomorphisms $u=h^{-1}\partial_{s}h$\,,  $v=h^{-1}\partial_{t}h$ and we put
\begin{equation}
{\cal R}= d''(h^{-1}d'h) + [\phi,\bar\phi_{h}]\,
\end{equation}
and $\Psi = i\,{\rm tr}[h^{-1}\tilde d h {\cal R}]$\,, where $\tilde d h = \partial_{s}h\,ds + \partial_{t}h\, dt$ 
is considered as the exterior differential of $h$ in the domain $\Delta$. It is convenient to rewrite $\Psi$ in the form
\begin{equation}
\Psi = i\,{\rm tr}[(u\,ds + v\,dt) {\cal R}]\,.
\end{equation}
Applying the Stokes formula to $\Psi$ (which is considered here as a 1-form in the domain $\Delta$) we get
\begin{equation}
\int_{\Delta}\tilde d \Psi = \int_{\partial \Delta} \Psi\,.                     \label{Stokes formula}
\end{equation}
The right hand side of the above expression can be computed straightforward from definition. In fact, after a short computation we obtain
\begin{equation}
\int_{\partial \Delta} \Psi = i \int_{a}^{b} {\rm tr} (v_{t}\cdot {\mathcal R}_{t}^{1,1})\,dt  +  Q_{2}^{1,1}(h_{a},k) - Q_{2}^{1,1}(h_{b},k)\,.
\end{equation}
Therefore, we need to show that the left hand side of (\ref{Stokes formula}) is an element in $d'A^{0,1} + d''A^{1,0}$, and hence, it suffices 
to show that $\tilde d  \Psi \in d'A^{0,1} + d''A^{1,0}$\,. \\

Now, from the definition of $\Psi$ we have 
\begin{eqnarray*}
\tilde d \Psi &=& i\, {\rm tr}[\tilde d (u\,ds + v\,dt)\,{\cal R} - (u\,ds + v\,dt) \tilde d {\cal R}] \\ 
              &=& i\, {\rm tr}[(\partial_{s}v - \partial_{t}u) {\cal R} - u\,\partial_{t}{\cal R} + v\,\partial_{s}{\cal R}]\,ds \wedge dt\,.
\end{eqnarray*}
On the other hand, a simple computation shows that 
\begin{equation}
\partial_{t}u = -vu + h^{-1}\partial_{t}\partial_{s}h\,, \quad\quad  \partial_{s}v = -uv + h^{-1}\partial_{s}\partial_{t}h\,, 
\end{equation}
\begin{equation}
\partial_{t}{\cal R} = d''D'v + [\phi,\partial_{t}\bar\phi_{h}]\,, \quad\quad  \partial_{s}{\cal R}= d''D'u + [\phi,\partial_{s}\bar\phi_{h}]\,. 
\end{equation}
Replacing these expressions in the formula for $\tilde d \Psi$ and writing ${\cal R}= R + [\phi,\bar\phi_{h}]$ we obtain 
\begin{eqnarray*}
\tilde d \Psi &=& i\, {\rm tr}\left[(vu-uv)R - u\,d''D'v + v\,d''D'u \right] ds\wedge dt \\ 
              && + \,i\, {\rm tr}\left[v\,[\phi,\partial_{s}\bar\phi_{h}] - 
              u\,[\phi,\partial_{t}\bar\phi_{h}]  + (vu-uv)[\phi,\bar\phi_{h}] \right] ds \wedge dt \,. 
\end{eqnarray*}
The first trace in the expression above does not depend on Higgs field $\phi$ (in fact, it is the same expression that is found in \cite{Kobayashi} for the ordinary case). 
The second trace is identically zero. In order to prove this, we need first to find explicit expressions for $\partial_{t}\bar\phi_{h}$ and $\partial_{s}\bar\phi_{h}$. 
Now (omitting the parameter $t$ for simplicity) we know from \cite{Simpson} that
\begin{equation}
\bar\phi_{h_{s+\delta s}} = u_{0}^{-1}\bar\phi_{h_{s}}u_{0} = \bar\phi_{h_{s}} + u_{0}^{-1}[\bar\phi_{h_{s}},u_{0}]
\end{equation}
where $u_{0}$ is a selfadjoint endomorphism such that $h_{s+\delta s}=h_{s}u_{0}$. Now
\begin{equation}
 h_{s+\delta s} = h_{s} + \partial_{s}h_{s}\cdot \delta s + {\cal O}(\delta s^{2})
\end{equation}
and hence, at first order in $\delta s$, we obtain $u_{0} = 1 + u\cdot\delta s$ and consequently $\partial_{s}\bar\phi_{h}=[\bar\phi_{h},u]$. In a similar way 
we obtain the formula $\partial_{t}\bar\phi_{h}=[\bar\phi_{h},v]$. Therefore, using these relations, the Jacobi identity and the cyclic property of the trace, 
we see that the second trace is identically zero. On the other hand, the term involving the curvature $R$ can be rewritten in terms of $u,v$ and their covariant derivatives. 
So, finally we get 
\begin{equation}
 \tilde d \Psi = - i\, {\rm tr}[v\, D'd''u + u\,d''D'v]\, ds\wedge dt\,.
\end{equation}
As it is shown in \cite{Kobayashi}, defining the (0,1)-form $\alpha = i\,{\rm tr}[v\,d''u]$ we obtain
\begin{equation}
 \tilde d \Psi = - [d'\alpha + d''\bar\alpha + i\,d''d' {\rm tr}(vu)]\, ds\wedge dt\,
\end{equation}    
and hence $\tilde d\Psi$ is an element in $d'A^{0,1} + d''A^{1,0}$\,. \;\;Q.E.D. \\

As a consequence of Lemma \ref{Main lemma L} we have an important result for piecewise differentiable closed curves. Namely, we have 

\begin{pro}\label{Main prop L} Let $h_{t}, \,\alpha \le t \le \beta\,,$ be a piecewise differentiable closed curve in 
{\rm Herm}$^{+}({\mathfrak E})$. Then
\begin{equation}
i \int_{\alpha}^{\beta} {\rm tr} \left(v_{t}\cdot {\mathcal R}_{t}^{1,1}\right)\,dt = 0 \quad   {\rm mod\,\,} d'A^{0,1} + d''A^{1,0}\,.
\end{equation} 
\end{pro}

\noindent{\it Proof:} Let $\alpha=a_{0}<a_{1} \cdots <a_{p}=\beta$ be the values of $t$ where $h_{t}$ is not differentiable. Now 
take a fixed point $k$ in ${\rm Herm}^{+}({\mathfrak E})$. Then, Lemma \ref{Main lemma L} applies for each triple $k,h_{a_{j}}, h_{a_{j+1}}$
with $j=0,1,...,p-1$ and the result follows.  \;\;Q.E.D.

\begin{cor}\label{} The Donaldson functional ${\mathcal L}(h,h')$ does not depend on the curve joining $h$ and $h'$\,. 
\end{cor}

\noindent{\it Proof:} Clearly, from the definition of $Q_{1}$
\begin{equation}
Q_{1}(h,h') + Q_{1}(h',h) = 0\,. 
\end{equation}
If $\gamma_{1}$ and $\gamma_{2}$ are two differentiable curves from $h$ to $h'$ and we apply Proposition \ref{Main prop L} 
to  $\gamma_{1}-\gamma_{2}$, we obtain
\begin{equation}
Q_{2}^{1,1}(h,h') + Q_{2}^{1,1}(h',h) = 0  \quad   {\rm mod\,\,} d'A^{0,1} + d''A^{1,0}\,,  
\end{equation}
and the result follows by integrating over $X$. \;\;Q.E.D.

\begin{pro}\label{} For any $h$ in ${\rm Herm}^{+}({\mathfrak E})$ and any constant $a>0$, the Donaldson functional 
satisfies ${\cal L}(h,ah)=0$\,.
\end{pro}

\noindent{\it Proof:} Clearly
\begin{equation}
Q_{1}(h,ah)={\rm log}\,{\rm det}[(ah)^{-1}h]=-r\,{\rm log\,}a\,. \nonumber
\end{equation}
Now, let $b={\rm log\,}a$ and consider the curve $h_{t}=e^{b(1-t)}h$ from $ah$ to $h$. For this curve $v_{t}=-bI$ and we have
\begin{equation}
{\cal R}_{t}^{1,1}=d''(h_{t}^{-1}d'h_{t}) + [\phi,\bar\phi_{t}] = d''(h^{-1}d'h) + [\phi,\bar\phi_{t}]\,, \nonumber
\end{equation} 
where $\bar\phi_{t}$ is an abbreviation for $\bar\phi_{h_{t}}$. Therefore, the (1,1)-component of $Q_{2}(h,ah)$ becomes 
\begin{equation}
Q_{2}^{1,1}(h,ah) = i\int_{0}^{1}{\rm tr}(v_{t}\cdot{\cal R}_{t}^{1,1})\,dt =i \int_{0}^{1}{\rm tr}\left[-b(R + [\phi,\bar\phi_{t}])\right] dt=-ib\,{\rm tr\,}R \nonumber
\end{equation} 
and hence, from the above we obtain
\begin{equation}
\frac{c}{n}\int_{X}Q_{1}(h,ah)\,\omega\wedge\omega^{n-1}/(n-1)! = -crb\,{\rm vol\,}X\,, \nonumber
\end{equation}
\begin{eqnarray*}
\int_{X} Q_{2}(h,ah)\wedge\omega^{n-1}/(n-1)! &=& \frac{-ib}{(n-1)!}\int_{X}{\rm tr}\,R\wedge \omega^{n-1} \\
&=&\frac{-2\pi b}{(n-1)!}\,{\rm deg\,}{\mathfrak E}
\end{eqnarray*}
and the result follows from the definition of the constant $c$\,.  \;\;Q.E.D.


\begin{lem}\label{Prop's of L} For any differentiable curve $h_{t}$ and any fixed point $k$ in {\rm Herm}$^{+}({\mathfrak E})$ 
we have
\begin{equation}
\partial_{t} Q_{1}(h_{t},k) = {\rm tr}(v_{t})\,,  
\end{equation}
\begin{equation}
\partial_{t} Q_{2}^{1,1}(h_{t},k) = i\,{\rm tr}(v_{t}\cdot{\mathcal R}_{t}^{1,1})
  \quad   {\rm mod\,\,} d'A^{0,1} + d''A^{1,0}\,.  
\end{equation}
\end{lem}

\noindent{\it Proof:} Since $k$ does not depend on $t$, we get
\begin{equation}
\partial_{t} Q_{1}(h_{t},k) =  \partial_{t}{\rm log}({\rm det}\, k^{-1}) + \partial_{t}{\rm log}({\rm det}\, h_{t}) = 
\partial_{t}{\rm log}({\rm det}\, h_{t}) = {\rm tr}(v_{t})\,. \nonumber
\end{equation}
Considering $b$ in (\ref{first Lemma}) as a variable, and differentiating that expression with respect to $b$, we 
obtain the formula.  \;\;Q.E.D. \\

By using the above Lemma, we have a formula for the derivative with respect to $t$ of Donaldson's functional
\begin{eqnarray*}
\frac{d}{dt}{\cal L}(h_{t},k) &=& \int_{X}\left[ i\,{\rm tr}(v_{t}\cdot{\cal R}_{t}^{1,1}) - \frac{c}{n}{\rm tr}(v_{t})\,
\omega\right]\wedge \frac{\omega^{n-1}}{(n-1)!}  \\
&=&   \int_{X}\left[{\rm tr}(v_{t}\cdot{\cal K}_{t}) - c\,{\rm tr}(v_{t})\right]\frac{\omega^{n}}{n!} \\
&=&   \int_{X}{\rm tr}\left[({\cal K}_{t} - cI)v_{t}\right]\frac{\omega^{n}}{n!} \,.
\end{eqnarray*}
Since $v_{t}=h_{t}^{-1}\partial_{t}h_{t}$ and we can consider the endomorphism ${\cal K}_{t}$ as an Hermitian form by defining ${\cal K}_{t}(s,s') = h_{t}(s,{\cal K}_{t}s')$, for any fixed Hermitian metric $k$ and any differentiable curve $h_{t}$ in 
Herm$^{+}({\mathfrak E})$ we obtain\footnote{Notice that from the 
definition (\ref{K endomorphism vs. form}), the endomorphism ${\cal K}_{t}$ can be written formally as ${\cal K}_{t}=h_{t}^{-1}{\cal K}_{t}(\cdot,\cdot)$ 
where ${\cal K}_{t}(\cdot,\cdot)$ denotes this time the mean curvature as a form. Therefore, we can express the derivative of the functional as an inner product 
of the forms ${\cal K}_{t} - c\,h_{t}$ and $\partial_{t}h_{t}$ as in (\ref{inner product}).} 
\begin{equation}
\frac{d}{dt}{\cal L}(h_{t},k) = \left({\cal K}_{t} - c\,h_{t},\partial_{t}h_{t}\right)\,,     \label{master eq.}
\end{equation}
where ${\cal K}_{t}$ is considered here as a form. For each $t$\,, we can consider $\partial_{t}h_{t}\in{\rm Herm}({\mathfrak E})$ as a tangent vector of 
${\rm Herm}^{+}({\mathfrak E})$ at $h_{t}$. Therefore, the differential $d{\cal L}$ of the functional evaluated at $\partial_{t}h_{t}$ is given by
\begin{equation}
d{\cal L}(\partial_{t}h_{t}) = \frac{d}{dt}{\cal L}(h_{t},k)\,,
\end{equation} 
and hence, the gradient of ${\cal L}$ (i.e., the vector field on ${\rm Herm}^{+}({\mathfrak E})$ dual to the form $d{\cal L}$ with respect to the invariant Riemannian 
metric introduced before) is given by $\nabla {\cal L} = {\cal K}-ch$\,. From the above analysis we conclude the following

\begin{thm}\label{} Let $k$ be a fixed element in ${\rm Herm}^{+}({\mathfrak E})$\,. Then, $h$ is a critical point of ${\cal L}$ if and only if ${\cal K} - c\,h=0$\,, 
i.e., if and only if $h$ is an Hermitian-Yang-Mills structure for ${\mathfrak E}$\,. 
\end{thm}
 
In order to derive some properties of ${\cal L}$ it is convenient to divide the Hichin-Simpson connection (see \cite{Simpson}, \cite{Simpson 2}) 
in the form ${\cal D}'_{h}=D'_{h} + \bar\phi_{h}$ and ${\cal D}''= D'' + \phi$\,. In fact, using the above decomposition it is not difficult to show that all 
critical points of ${\cal L}$ correspond to an absolute minimum. 

\begin{thm}\label{Minimum of L} Let $k$ be a fixed Hermitian structure of ${\mathfrak E}$ and $\tilde h$ a critical point of ${\cal L}(h,k)$\,, 
then the Donaldson functional attains an absolute minimum at $\tilde h$\,.
\end{thm}

\noindent{\it Proof:} Let $h_{t}\,, 0\le t\le 1\,,$ be a differentiable curve such that $h_{0}=\tilde h$, then we can compute 
straightforward the second derivative of ${\cal L}$ 
\begin{eqnarray*}
\frac{d^{2}}{dt^{2}}{\cal L}(h_{t},k) &=& \frac{d}{dt}\int_{X}{\rm tr}\left[({\cal K}_{t} - cI)v_{t}\right]\frac{\omega^{n}}{n!} \\
 &=&\int_{X}{\rm tr}\left[\partial_{t}{\cal K}_{t}\cdot v_{t} + ({\cal K}_{t} - cI)\partial_{t}v_{t}\right]\frac{\omega^{n}}{n!}\,.
\end{eqnarray*}
Since $h_{0}$ is a critical point of the functional, then ${\cal K}_{t} - cI = 0$ at $t=0$\,, and hence
\begin{equation}
\frac{d^{2}}{dt^{2}}{\cal L}(h_{t},k)\big|_{t=0} = \int_{X}{\rm tr}(\partial_{t}{\cal K}_{t}\cdot v_{t})\,\frac{\omega^{n}}{n!} \big|_{t=0}\,.
\end{equation}
On the other hand, $\partial_{t}{\cal K}_{t}$ can be written in terms of the endomorphism $v_{t}$ in the following way
\begin{eqnarray*}
{\cal D''}{\cal D'}v_{t}&=&{\cal D''}(D'v_{t} + [\bar\phi_{t},v_{t}])\\
&=& D''D'v_{t} + [\phi,D'v_{t}] + D''[\bar\phi_{t},v_{t}] + [\phi,[\bar\phi_{t},v_{t}]]\,,
\end{eqnarray*}
and since $\partial_{t}\phi_{t}=[\bar\phi_{t},v_{t}]$ we get 
\begin{equation}
\partial_{t}{\cal R}_{t}^{1,1} = \partial_{t}R_{t} + [\phi,\partial_{t}\bar\phi_{t}] = D''D'v_{t} + [\phi,[\bar\phi_{t},v_{t}]]\,.
\end{equation}
Therefore, taking the trace with respect to $\omega$ (i.e., applying the $i\Lambda$ operator) we obtain
\begin{equation}
i\Lambda{\cal D''}{\cal D'}v_{t}=i\Lambda\partial_{t}{\cal R}_{t}^{1,1}=\partial_{t}{\cal K}_{t}\,.   \label{box op.} 
\end{equation}
Hence, replacing this in the expression for the second derivative of ${\cal L}$ we find
\begin{equation}
\frac{d^{2}}{dt^{2}}{\cal L}(h_{t},k)\big|_{t=0} = \int_{X}{\rm tr}(i\Lambda{\cal D''}{\cal D'}v_{t}\cdot v_{t})\frac{\omega^{n}}{n!} \big|_{t=0} = \lVert{\cal D'}v_{t}\rVert^{2}_{t=0}\,,
\end{equation}
(that is, $h_{0}$ must be at least a local minimum of ${\cal L}$). Now suppose in addition that $h_{1}$ is an arbitrary element in 
${\rm Herm}^{+}({\mathfrak E})$ and joint them by a geodesic $h_{t}$, and hence $\partial_{t}v_{t}=0$\,. Therefore, for a such a geodesic we have 
\begin{equation}
\frac{d^{2}}{dt^{2}}{\cal L}(h_{t},k) = \int_{X}{\rm tr}(\partial_{t}{\cal K}_{t}\cdot v_{t})\,\frac{\omega^{n}}{n!}\,.
\end{equation}
Following the same procedure we have done before, but this time at $t$ arbitrary, we get for $0\le t\le 1$
\begin{equation}
\frac{d^{2}}{dt^{2}}{\cal L}(h_{t},k) = \lVert{\cal D'}v_{t}\rVert^{2}_{h_{t}} \ge 0
\end{equation}
(since there is an implicit dependence on $t$ on the right hand side via ${\cal D'}$, we write a subscript $h_{t}$ in the norm) and 
it follows that ${\cal L}(h_{0},k)\le{\cal L}(h_{1},k)$. Now if we assume that $h_{1}$ is also a critical point of ${\cal L}$\,, 
we necessarily obtain the equality. Therefore, it follows that the minimum defined for any critical point of ${\cal L}$ is an absolute minimum. \;\;Q.E.D. \\


Let $k$ be a fixed Hermitian structure, then any Hermitian metric $h$ will be of the form $ke^{v}$ for some section $v$ of ${\rm End}(E)$ over $X$. 
We can join $k$ to $h$ by the geodesic $h_{t}=ke^{tv}$ where $0\le t\le 1$ (note that here $v_{t}=h_{t}^{-1}\partial_{t}h_{t}=v$ 
is constant, i.e., it does not depend on $t$). Now, in the proof of Theorem \ref{Minimum of L} we got an expression for the second derivative 
${\cal L}(h_{t},k)$ for any curve $h_{t}$. Namely 
\begin{equation}
\frac{d^{2}}{dt^{2}}{\cal L}(h_{t},k) = \int_{X}{\rm tr}\left[\partial_{t}{\cal K}_{t}\cdot v_{t} + ({\cal K}_{t} - cI)\partial_{t}v_{t}\right]\frac{\omega^{n}}{n!}\,.
\end{equation}
Notice that in our case, the chosen curve is such that $h_{0}=k$, since it is also a geodesic $\partial_{t}v_{t}=0$ we have  
\begin{equation}
\frac{d^{2}}{dt^{2}}{\cal L}(h_{t},k)=\int_{X}{\rm tr}(\partial_{t}{\cal K}_{t}\cdot v)\frac{\omega^{n}}{n!}=\lVert{\cal D'}v\rVert^{2}_{h_{t}}\,. \label{2nd derivative of L}
\end{equation}
Therefore, following \cite{Siu}, the idea is to find a simple expression for $\lVert{\cal D'}v\rVert^{2}_{h_{t}}$ or equivalently for  $\lVert{\cal D''}v\rVert^{2}_{h_{t}}$ 
and to integrate it twice with respect to $t$\,. We can do this using local coordinates, in fact, at any point in $X$ we can choose a local frame field 
so that $h_{0}=I$ and $v={\rm diag}(\lambda_{1},...,\lambda_{r})$\,. In particular, using such a local frame field we have $h_{t}^{ij}=e^{-\lambda_{j}t}\delta_{ij}$\,,
and hence (after a short computation) we obtain
\begin{equation}
\lVert{\cal D''}v\rVert^{2}_{h_{t}}=\int_{X}\sum_{i,j=1}^{r}e^{(\lambda_{i}-\lambda_{j})t}\, |{\cal D''}v_{j}^{i}|^{2}\, \frac{\omega^{n-1}}{(n-1)!} \,. 
\end{equation}
Now, at $t=0$ the functional ${\cal L}(h_{t},k)$ vanishes and since $k=h_{0}$ is not necessarily an 
Hermitian-Yang-Mills structure, we have
\begin{equation}
\frac{d}{dt}{\cal L}(h_{t},k)\big|_{t=0} = \int_{X}{\rm tr}\left[({\cal K}_{0} - cI)v\right]\frac{\omega^{n}}{n!} \,.
\end{equation}
Then, by integrating twice (\ref{2nd derivative of L}) we obtain
\begin{equation}
{\cal L}(h_{t},k)= t\int_{X}{\rm tr}\left[({\cal K}_{0} - cI)v\right]\frac{\omega^{n}}{n!} + \int_{X}\sum_{i,j=1}^{r}\psi_{t}(\lambda_{i},\lambda_{j})|{\cal D''}v_{j}^{i}|^{2}\, \frac{\omega^{n-1}}{(n-1)!} \label{L Simpson}
\end{equation}
where $\psi_{t}$ is a function given by 
\begin{equation}
\psi_{t}(\lambda_{i},\lambda_{j})=\frac{e^{(\lambda_{i}-\lambda_{j})t} - (\lambda_{i}-\lambda_{j})t - 1}
{(\lambda_{i}-\lambda_{j})^{2}}\,.
\end{equation} 
In particular, at $t=1$ the expression (\ref{L Simpson}) corresponds (up to a constant term) to the definition of Donaldson's functional 
given by Simpson in \cite{Simpson}. Notice also that if the initial metric $k=h_{0}$ is Hermitian-Yang-Mills, the first term of the right hand side of 
(\ref{L Simpson}) vanishes and the functional coincides with the Donaldson functional used by Siu in \cite{Siu}.

\section{The evolution equation}
For the construction of Hermitian-Yang-Mills structures, the standard procedure is to start with a fixed Hermitian metric $h_{0}$ and try 
to find from it an Hermitian metric satisfying ${\cal K}=cI$ using a curve $h_{t}$\,, $0\le t<\infty$\, (in other words, we try to find that metric by deforming $h_{0}$ 
through 1-parameter family of Hermitian metrics) and we expect that at $t=\infty$, the metric will be Hermitian-Yang-Mills. \\

At this point, it is convenient to introduce the operator $\tilde\Box_{h} = i\Lambda{\cal D''}{\cal D'}_{h}$, which depends on the metric $h$\,. 
Using it, we can rewrite (\ref{box op.}) as $\partial_{t}{\cal K}_{t}=\tilde\Box_{t}v_{t}\,,$ where the subscript $t$ reminds the dependence 
of the operator on the metric $h_{t}$\,. \\

As we said before, to get an Hermitian-Yang-Mills metric we want to make ${\cal K} -cI$ vanish. Therefore, a natural choice is to go along the 
global gradient direction of the functional given by the global $L^{2}$-norm of ${\cal K}_{t} -cI$\,. Therefore, taking the derivative of this 
functional we obtain
\begin{eqnarray*}
\frac{d}{dt}\lVert{\cal K}_{t} -cI\rVert^{2}  &=& \int_{X} 2\, {\rm tr}\left(\partial_{t}{\cal K}_{t}\cdot({\cal K}_{t} -cI)\right) \frac{\omega^{n}}{n!}\\
 &=& 2 \int_{X} {\rm tr}\left(\tilde\Box_{t}v_{t}\cdot({\cal K}_{t} -cI)\right) \frac{\omega^{n}}{n!} \\
 &=& 2 \int_{X} {\rm tr}\left(v_{t}\cdot\tilde\Box_{t}{\cal K}_{t}\right) \frac{\omega^{n}}{n!}\,, 
\end{eqnarray*}
and the equation that naturally emerges (i.e., the associated steepest descent curve) is $v_{t}=-\tilde\Box_{t}{\cal K}_{t}$\,, or equivalently
\begin{equation}
h_{t}^{-1}\partial_{t}h_{t} = -i\Lambda{\cal D''}{\cal D'}_{h}{\cal K}_{t}\,.
\end{equation}
Since ${\cal K}_{t}$ is of degree two, the right hand side of the above equation becomes a term of degree four and hence we get 
at the end a nonlinear equation of degree four. To do the analysis, it is easier to deal with an equation of lower degree. In fact, this is one of the 
reasons for introducing the Donaldson functional. Following the same argument we did before, but this time using the functional ${\cal L}(h_{t},k)$ with $k$ fixed, 
in place of the functional $\lVert{\cal K}_{t} -cI\rVert^{2}$\,, we end up with a nonlinear equation of degree two (the heat equation),
to be more precise, we obtain directly from (\ref{master eq.}) the equation 
\begin{equation}
\partial_{t}h_{t} = -( {\cal K}_{t}-ch_{t})\,,  
\end{equation}
where this time, ${\cal K}_{t}$ represents the associated two form, and not an endomorphism\footnote{Notice that the equivalent equation involving 
endomorphisms will be $v_{t} = -( {\cal K}_{t}-cI)\,.$}. Simpson has shown that also for the Higgs case, 
we have always solutions of the above non linear evolution equation. This was proved in \cite{Simpson} for the non-compact case satisfying some additional 
conditions. That proof covers the compact K\"ahler case without any change. Then, from \cite{Simpson} we have the following

\begin{thm}\label{} Given an Hermitian structure $h_{0}$ on ${\mathfrak E}$\,, the non-linear evolution equation    
\begin{equation}
\partial_{t}h_{t} = -( {\cal K}_{t}-ch_{t})\,
\end{equation}
has a unique smooth solution defined for $0\le t < \infty$\,. 
\end{thm}

In the rest of this section, we study some properties of the solutions of the evolution equation. In particular, we are interested in the study of the 
mean curvature when the paramater $t$ goes to infinity.  

\begin{pro}\label{} Let $h_{t}\,, 0\le t < \infty\,,$ be a 1-parameter family of ${\rm Herm}^{+}({\mathfrak E})$ satisfying the evolution equation. 
Then \\
\noindent {\bf (i)} For any fixed Hermitian structure $k$ of ${\mathfrak E}$\,, the functional ${\cal L}(h_{t},k)$ is a monotone decreasing function 
of $t$; that is 
\begin{equation}
\frac{d}{dt} {\cal L}(h_{t},k) =  -\lVert{\cal K}_{t}-cI\rVert^{2} \le 0 \,;
\end{equation}
{\bf (ii)} $\max |{\cal K}_{t}-cI|^{2}$ is a monotone decreasing function of $t$\,;\\
{\bf (iii)} If ${\cal L}(h_{t},k)$ is bounded below, i.e., ${\cal L}(h_{t},k)\ge A > -\infty$ for some real constant $A$ and $0\le t < \infty$\,, 
then
\begin{equation}
\max_{X}|{\cal K}_{t}-cI|^{2} \rightarrow 0 \quad as \quad t\rightarrow \infty\,.                      \label{item iii} 
\end{equation}
\end{pro}

\noindent{\it Proof:} From the proof of Lemma \ref{Prop's of L} we know that 
\begin{equation}
\frac{d}{dt}{\cal L}(h_{t},k) = \left({\cal K}_{t} - c\,h_{t},\partial_{t}h_{t}\right)\,.
\end{equation}
Since $h_{t}$ is a solution of the evolution equation, we get 
\begin{equation}
\frac{d}{dt}{\cal L}(h_{t},k) = -\left({\cal K}_{t} - c\,h_{t},{\cal K}_{t} - c\,h_{t}\right)= - \lVert{\cal K}_{t} - c\,h_{t}\rVert^{2}
\end{equation}
and (i) follows from the definition of the Riemannian structure in ${\rm Herm}^{+}({\mathfrak E})$ (considered this time as a metric for endomorphisms). The proofs 
of (ii) and (iii) are similar to the proof in the classical case \cite{Kobayashi}, but we need to work this time with the 
operator $\tilde\Box_{h} = i\Lambda{\cal D}''{\cal D}'_{h}$ instead of the operator $\Box_{h} = i\Lambda D''D'_{h}$\,. In fact, from this 
definition $\tilde\Box v_{t}=\partial_{t}{\cal K}_{t}$ and since $v_{t}=-({\cal K}_{t}-cI),$ we get 
\begin{equation}
(\tilde\Box + \partial _{t}){\cal K}_{t}=0\,. 
\end{equation}
On the other hand,
\begin{eqnarray*}
{\cal D''}{\cal D'}|{\cal K}_{t}-cI|^{2} &=& {\cal D''}{\cal D'} {\rm tr}({\cal K}_{t}-cI)^{2}\\
&=& 2\,{\rm tr}(({\cal K}_{t}-cI){\cal D''}{\cal D'}{\cal K}_{t}) + 2\,{\rm tr}({\cal D''}{\cal K}_{t}\cdot{\cal D'}{\cal K}_{t})\,
\end{eqnarray*}
and by taking the trace with respect to $\omega$ we get
\begin{eqnarray*}
\tilde\Box|{\cal K}_{t}-cI|^{2} &=& 2\,{\rm tr}(({\cal K}_{t}-cI)\tilde\Box{\cal K}_{t}) + 2\,i\Lambda\, {\rm tr}({\cal D''}{\cal K}_{t}\cdot{\cal D'}{\cal K}_{t}) \\
&=& - 2\,{\rm tr}(({\cal K}_{t}-cI)\partial_{t}{\cal K}_{t}) - 2\,|{\cal D''}{\cal K}_{t}|^{2}\\
&=& - \partial_{t}|{\cal K}_{t}-cI|^{2} - 2\,|{\cal D''}{\cal K}_{t}|^{2}\,.
\end{eqnarray*}
So, finally we obtain
\begin{equation}
(\partial_{t} + \tilde\Box)|{\cal K}_{t}-cI|^{2} = - 2\,|{\cal D''}{\cal K}_{t}|^{2} \le 0             \label{function f} 
\end{equation} 
and (ii) follows from the maximum principle. Finally, (iii) follows from (ii) and (i) in a similar way to the classical case 
(see \cite{Kobayashi} for details). \;\;Q.E.D. \\
 
At this point we introduce the main result of the section. This establishes a relation among the boundedness property of 
Donaldson's functional, the semistability and the existence of approximate Hermitian-Yang-Mills structures.

\begin{thm}\label{} Let $\mathfrak E$ be a Higgs bundle over a compact K\"ahler manifold $X$ with K\"ahler form $\omega$\,. Then we have the 
implications ${\rm (i)} \rightarrow {\rm (ii)} \rightarrow {\rm (iii)}$ for the following statements: \\
\noindent {\bf (i)} for any fixed Hermitian structure $k$ in ${\mathfrak E}$\,, there exists a constant $B$ such that ${\cal L}(h,k)\ge B$ for all 
Hermitian structures $h$ in ${\mathfrak E}$\,; \\
\noindent {\bf (ii)} ${\mathfrak E}$ admits an approximate Hermitian-Yang-Mills structure, i.e., given $\epsilon > 0$ there exists an Hermitian 
structure $h$ in ${\mathfrak E}$ such that
\begin{equation}
 \max_{X}|{\cal K} -ch| < \epsilon\,;
\end{equation}
\noindent {\bf (iii)} ${\mathfrak E}$ is $\omega$-semistable\,.
\end{thm}

\noindent{\it Proof:} Assume (i). Then the funcional is bounded below by a constant $A$. In particular ${\cal L}(h_{t},k)\ge A$ for $h_{t}$\,, $0\le t < \infty$\,, 
a one-parameter family of Hermitian structures satisfying the evolution equation. Therefore, from (\ref{item iii}) it follows that given $\epsilon >0$\, there exists $t_{0}$ such that\footnote{Notice that 
given $\epsilon >0$, any metric $h=h_{t_{1}}$ with $t_{1}>t_{0}$ in principle satisfies (\ref{ineq.}).} 
\begin{equation}
\max_{X}|{\cal K}_{t}-cI| < \epsilon \quad {\rm for} \quad t>t_{0}\,.   \label{ineq.}
\end{equation}
This shows that (i) implies (ii). On the other hand, (ii) $\rightarrow$ (iii) has been proved by Bruzzo and Gra\~na-Otero in 
\cite{Bruzzo-Granha}. \;\;Q.E.D.

\section{Semistable Higgs bundles}
We need some results which allow us to solve some problems about Higgs bundles by induction on the rank. This section is 
essentially a natural extension to Higgs bundles of the classical case, which is explained in detail in \cite{Kobayashi}. \\     

Let
\begin{equation}
 \xymatrix{
0 \ar[r]  &  {\mathfrak E'} \ar[r]^{\iota}  &   {\mathfrak E} \ar[r]^{p}  &   {\mathfrak E''} \ar[r] &  0  
}   \label{model exact seq}
\end{equation}
be an exact sequence of Higgs bundles over a K\"ahler manifold. As in the ordinary case, an Hermitian structure $h$ in ${\mathfrak E}$ 
induces Hermitian structures $h'$ and $h''$ in  ${\mathfrak E'}$  and ${\mathfrak E''}$ respectively. We have also a second fundamental form 
$A_{h}\in A^{1,0}({\rm Hom}(E',E''))$ and its adjoint $B_{h}\in A^{0,1}({\rm Hom}(E'',E'))$ where, as usual, $B_{h}^{*}=-A_{h}$\,. In a similar way, 
some properties which hold in the ordinary case, also hold in the Higgs case. \\
 
\begin{pro}\label{Fundamental properties Q1,Q2} Given an exact sequence $(\ref{model exact seq})$ and a pair of Hermitian structures $h,k$ in ${\mathfrak E}$\,. 
Then the function $Q_{1}(h,k)$ and the form $Q_{2}(h,k)$ satisfies the following properties:\\
\noindent {\bf (i)} $Q_{1}(h,k) = Q_{1}(h',k') + Q_{1}(h'',k'')\,,$\\
\noindent{\bf (ii)} $Q_{2}(h,k) = Q_{2}(h',k') + Q_{2}(h'',k'') - i\,{\rm tr}[B_{h}\wedge B_{h}^{*} - B_{k}\wedge B_{k}^{*}]$\\
$\;\; {\rm mod}\, d'A^{0,1} + d''A^{1,0}$\,.\\
\end{pro}

\noindent {\it Proof:} (i) is straightforward from the definition of $Q_{1}$\,. On the other hand, (ii) follows from an analysis similar 
to the ordinary case. \\

Since the sequence (\ref{model exact seq}) is in particular an exact sequence of holomorphic vector bundles, for any $h$ we have a splitting 
of the exact sequence by $C^{\infty}$-homomorphisms $\mu_{h}: E\rightarrow E'$ and $\lambda_{h}: E''\rightarrow E$\,. In particular, 
\begin{equation}
B_{h}=\mu_{h}\circ d''\circ \lambda_{h}\,.   
\end{equation}

We consider now a curve of Hermitian structures $h=h_{t}$\,, $0\le t \le 1$ such that $h_{0}=k$ and $h_{1}=h$\,. Corresponding to $h_{t}$ 
we have a family of homomorphisms $\mu_{t}$ and $\lambda_{t}$\,. We define the homomorphism $S_{t}: E''\rightarrow E'$ given by
\begin{equation}
\lambda_{t}-\lambda_{0}=\iota\circ S_{t}\,. 
\end{equation}
A short computation (see \cite{Kobayashi}, Ch.VI) shows that $\partial_{t}B_{t}=d''(\partial_{t}S_{t})$ and choosing convenient orthonormal local frames for ${\mathfrak E}'$ and 
${\mathfrak E}''$, we know that the endomorphism $v_{t}$ can be represented by the matrix
\[v_{t} = \left( \begin{array}{cc} v'_{t} & -\partial_{t}S_{t} \\
                                 -(\partial_{t}S_{t})^{*} & v''_{t} \\
               \end{array} \right)\]
where $v'_{t}, v''_{t}$ are the natural endomorphisms associated to $h'_{t}, h''_{t}$ respectively. Now, from the ordinary case we have  

\[R_{t} = \left( \begin{array}{cc} R'_{t} -B_{t}\wedge B_{t}^{*} & D'B_{t} \\
                                 -D''B_{t}^{*} &  R''_{t} -B_{t}^{*}\wedge B_{t} \\
               \end{array} \right),\]
where $R'_{t}$ and $R''_{t}$ are the Chern curvatures of ${\mathfrak E}'$ and ${\mathfrak E}''$ associated to the metrics 
$h'_{t}$ and $h''_{t}$ respectively. Now ${\cal R}^{1,1}_{t}=R_{t}+[\phi,\phi_{t}]$ and since ${\mathfrak E}'$ and ${\mathfrak E}''$ are Higgs 
subbundles of ${\mathfrak E}$, we obtain a simple expression for the (1,1)-component of the Hitchin-Simpson curvature

\[{\cal R}^{1,1}_{t} = \left( \begin{array}{cc} {\cal R}'^{1,1}_{t} -B_{t}\wedge B_{t}^{*}  & D'B_{t} \\
                                 -D''B_{t}^{*} &  {\cal R}''^{1,1}_{t} -B_{t}^{*}\wedge B_{t}  \\
               \end{array} \right),\]
where $ {\cal R}'^{1,1}_{t} = R'_{t} + [\phi,\phi_{t}]_{E'}$ and $ {\cal R}''^{1,1}_{t}= R''_{t} + [\phi,\phi_{t}]_{E''}$. Hence, at this point 
we can compute the trace
\begin{eqnarray*}
{\rm tr}(v_{t}\cdot{\cal R}^{1,1}_{t}) &=& {\rm tr}(v'_{t}\cdot{\cal R}'^{1,1}_{t}) +  {\rm tr}(v''_{t}\cdot{\cal R}''^{1,1}_{t}) \\
                                  && + {\rm tr}(\partial_{t}S_{t}\cdot D''B_{t}^{*}) -  {\rm tr}((\partial_{t}S_{t})^{*}\cdot D'B_{t})\\
                                  && + {\rm tr}(v'_{t}\cdot B_{t}\wedge B_{t}^{*}) - {\rm tr}(v''_{t}\cdot B^{*}_{t}\wedge B_{t})\,.
\end{eqnarray*}
The last four terms are exactly the same as in the ordinary case. Finally we get that, modulo an element in $d'A^{0,1}+d''A^{1,0}$
\begin{equation}
 {\rm tr}(v_{t}\cdot{\cal R}^{1,1}_{t}) = {\rm tr}(v'_{t}\cdot{\cal R}'^{1,1}_{t}) +  {\rm tr}(v''_{t}\cdot{\cal R}''^{1,1}_{t}) 
- \partial_{t} {\rm tr}(B_{t}\wedge B_{t}^{*})\,.
\end{equation}
Then, multiplying the last expression by $i$ and integrating from $t=0$ to $t=1$ we obtain (ii). \,\, Q.E.D.  \\

From Proposition \ref{Fundamental properties Q1,Q2} we get an important result for the compact case when 
$\mu({\mathfrak E})=\mu({\mathfrak E}')=\mu\,.$  Indeed, in that case we have also $\mu({\mathfrak E}'')=\mu\,.$ Then, by integrating 
$Q_{1}(h,k)$ and $Q_{2}(h,k)$ over $X$, and since 
\begin{equation}
-i{\rm tr}(B\wedge B^{*})\wedge \omega^{n-1} = |B|^{2}\omega^{n}/n!  
\end{equation}
we have the same constant $c$ for all functionals ${\cal L}(h,k)\,,{\cal L}(h',k')$ and ${\cal L}(h'',k'')$ and we obtain the following

\begin{cor}\label{} Given an exact sequence $(\ref{model exact seq})$ over a compact K\"ahler manifold $X$ with 
$\mu({\mathfrak E}) = \mu({\mathfrak E}')$ and a pair of Hermitian structures $h$ and $k$ in ${\mathfrak E}$, the functional 
${\cal L}(h,k)$ satisfies the following relation
\begin{equation}
 {\cal L}(h,k) = {\cal L}(h',k') + {\cal L}(h'',k'') + \lVert B_{h}\rVert^{2} - \lVert B_{k}\rVert^{2}\,.
\end{equation}
\end{cor}

In the one-dimensional case, when $X$ is a compact Riemann surface, the notion of stability (resp. semistability) does not 
depend on the K\"ahler form $\omega$, therefore we can establish our results without make reference to any $\omega$. At this point we can
establish a boundedness property for the Donaldson functional for semistable Higgs bundles over Riemann surfaces. To be precise we have

\begin{thm}\label{Boundedness of L} Let ${\mathfrak E}$ be a Higgs bundle over a compact Riemann surface $X$. If it is semistable, then for 
any fixed Hermitian structure $k$ in ${\rm Herm}^{+}(\mathfrak E)$ the set $\{{\cal L}(h,k), h\in{\rm Herm}^{+}(\mathfrak E)\}$ is bounded below.
\end{thm}

\noindent {\it Proof:} Fix $k$ and assume that ${\mathfrak E}$ is semistable. The proof runs by induction on the rank of ${\mathfrak E}$. 
If it is stable, then by \cite{Simpson} there exists an Hermitian-Yang-Mills structure $h_{0}$ on it, and we know that Donaldson's functional 
must attain an absolute minimum at $h_{0}\,,$ i.e., for any other metric $h$  
\begin{equation}
{\cal L}(h,k)\ge {\cal L}(h_{0},k) 
\end{equation}
and hence the set is bounded below. Now, suppose ${\mathfrak E}$ is not stable, then among all proper non-trivial Higgs subsheaves with torsion-free quotient and 
the same slope as ${\mathfrak E}$ we choose one, say ${\mathfrak E}'$, with minimal rank. Since $\mu({\mathfrak E}')=\mu({\mathfrak E})$ this sheaf 
is necessarily stable\footnote{If ${\mathfrak E}'$ is not stable, there exists a proper Higgs subsheaf 
${\mathfrak F}'$ of ${\mathfrak E}'$ with $\mu({\mathfrak F}') \ge \mu({\mathfrak E}')$ and since ${\mathfrak F}'$ is clearly a subsheaf of 
${\mathfrak E}$ and this is semistable we necessarily obtain $\mu({\mathfrak F}') = \mu({\mathfrak E})$, which is a contradiction, 
because ${\mathfrak E}'$ was chosen with minimal rank.}. Now let ${\mathfrak E}''= {\mathfrak E}/{\mathfrak E}'$, then using Lemma 7.3 in 
\cite{Kobayashi} it follows that $\mu({\mathfrak E}'') = \mu({\mathfrak E})$  and ${\mathfrak E}''$ is semistable\footnote{In fact, if ${\mathfrak E}''$ is not 
semistable, then there exists a proper Higgs subsheaf $\mathfrak H$ of ${\mathfrak E}''$ with 
$\mu(\mathfrak H)>\mu({\mathfrak E}'')$. Then, using Lemma 7.3 in \cite{Kobayashi} we have 
$\mu({\mathfrak E}'') > \mu({\mathfrak E}''/{\mathfrak H})$. Defining $\mathfrak K$ as the kernel of the morphism
${\mathfrak E}\longrightarrow {\mathfrak E}''/{\mathfrak H}$, we get the exact sequence 
\begin{displaymath}
 \xymatrix{
0   \ar[r]  &  {\mathfrak K} \ar[r]  &   {\mathfrak E} \ar[r]  &   {\mathfrak E}''/{\mathfrak H}  \ar[r]  &    0 
}
\end{displaymath}
and since $\mu(\mathfrak E)=\mu({\mathfrak E}'')$, using again the same Lemma in \cite{Kobayashi} we conclude that 
$\mu(\mathfrak K)>\mu(\mathfrak E)$, which contradicts the semistability of $\mathfrak E$.} and hence we have the following exact sequence of sheaves 
\begin{equation}
 \xymatrix{
0 \ar[r]  &  {\mathfrak E'} \ar[r]  &   {\mathfrak E} \ar[r]  &   {\mathfrak E''} \ar[r] &  0 
}   \label{model exact seq 2}
\end{equation}
where ${\mathfrak E}'$ and ${\mathfrak E}''$ are both torsion-free. Since ${\rm dim\,}X=1$ they are also locally free and hence the sequence is in fact an 
exact sequence of Higgs bundles. Assume now that $h$ is an arbitrary metric on ${\mathfrak E}$, then by applying the preceding Corollary to the metrics $h$ and $k$ we obtain 
\begin{equation}
 {\cal L}(h,k) = {\cal L}(h',k') + {\cal L}(h'',k'') + \lVert B_{h}\rVert^{2} - \lVert B_{k}\rVert^{2}\,,
\end{equation}
where $h',k'$ and $h'',k''$ are the Hermitian structures induced by $h,k$ in ${\mathfrak E}'$ and  ${\mathfrak E}''$ respectively. If the rank of 
${\mathfrak E}$ is one, it is stable and hence ${\cal L}(h,k)$ is bounded below by a constant which depends on $k$. Now, by the inductive hypothesis,  
${\cal L}(h',k')$ and ${\cal L}(h'',k'')$ are bounded below by constants depending only on $k'$ and $k''$ respectively. Then
${\cal L}(h,k)$ is bounded below by a constant depending on $k$.  \;\;Q.E.D. \\

As a consequence of the above result, we get that in the one-dimensional case all three conditions in the main theorem are equivalent. As a consequence 
we obtain the following 
\begin{cor}\label{Main Corollary} Let ${\mathfrak E}$ be a Higgs bundle over a compact Riemann surface $X$. 
Then ${\mathfrak E}$ is semistable if and only if ${\mathfrak E}$ admits an approximate Hermitian-Yang-Mills structure. 
\end{cor}
This equivalence between the notions of approximate Hermitian-Yang-Mills structures and semistability is one version of the so called 
{\it Hitchin-Kobayashi correspondence} for Higgs bundles. As a consequence of the Corollary \ref{Main Corollary} we see that in the one-dimensional case, 
all results about Higgs bundles written in terms of approximate Hermitian-Yang-Mills structures can be traslated in terms of semistability. In 
particular we have
\begin{cor}\label{semsitability property} If ${\mathfrak E}_{1}$ and ${\mathfrak E}_{2}$ are semistable Higgs bundles over a compact 
Riemann surface $X$, then so is the tensor product ${\mathfrak E}_{1}\otimes{\mathfrak E}_{2}$. Furthermore 
if $\mu({\mathfrak E}_{1})=\mu({\mathfrak E}_{2})$, so is the Whitney sum 
${\mathfrak E}_{1}\oplus{\mathfrak E}_{2}$\,.
\end{cor}
\begin{cor}\label{} If $\mathfrak E$ is semistable Higgs bundle over a Riemann surface $X$, then so is the tensor product bundle 
${\mathfrak E}^{\otimes p}\otimes{\mathfrak E}^{*\otimes q}$ and the exterior product bundle $\bigwedge^{p}{\mathfrak E}$ whenever $p\le r$\,. 
\end{cor}
The equivalence between the existence of approximate Hermitian-Yang-Mills structures and semistability is also true in higher dimensions. However, 
since torsion-free sheaves over compact K\"ahler manifolds with ${\rm dim\,}X\ge 2$ may not be locally free (they are locally free 
only outside its singularity set) it is necessary to use Higgs sheaves and not only Higgs bundles.

\section{Admissible metrics for Higgs sheaves}

A natural notion of a metric on a torsion-free sheaf is that of admissible Hermitian structure. This was first introduced by S. Bando and Y.-T. Siu 
in \cite{Bando-Siu}. In their article, they proved first the existence of admissible structures on any torsion-free 
sheaf, and then obtained an equivalence between the stability of a torsion-free sheaf and the existence of an admissible Hermitian-Yang-Mills metric 
on it, thus extending the Hitchin-Kobayashi correspondence to torsion-free sheaves. \\

Admissible structures were used again by I. Biswas and G. Schumacher \cite{Biswas-Schumacher} to prove an extended version of the correspondence 
of Bando and Siu to the Higgs case. In this last section, we briefly discuss some of these notions. \\   

Let ${\mathfrak E}$ be a torsion-free Higgs sheaf over a compact K\"ahler manifold $X$. The singularity set of $\mathfrak E$ is the subset 
$S=S(\mathfrak E)\subset X$ where $\mathfrak E$ is not locally free. As is well known, $S$ is a complex analytic subset with ${\rm codim}S\ge2$. Following \cite{Biswas-Schumacher}, \cite{Bando-Siu}, an 
{\it admissible structure} $h$ on $\mathfrak E$ is an Hermitian metric on the bundle ${\mathfrak E}|_{X\backslash S}$ with the following two 
properties: \\

\noindent{\bf (i)} The Chern curvature $R$ of $h$ is square-integrable, and \\
{\bf (ii)} The corresponding mean curvature $K=i\Lambda R$ is $L^{1}$-bounded. \\   

Let consider now the natural embedding of $\mathfrak E$ into its double dual ${\mathfrak E}^{\vee\vee}$; since 
$S({\mathfrak E}^{\vee\vee})\subset S(\mathfrak E)$, an admissible structure on ${\mathfrak E}^{\vee\vee}$ restrics to an admissible 
structure on $\mathfrak E$. An admissible structure $h$ is called an {\it admissible Hermitian-Yang-Mills structure} if on $X\backslash S$ 
the mean curvature of the Hitchin-Simpson connection is proportional to the identity. In other words if
\begin{equation}
{\mathcal K} = K + i\Lambda [\phi,\bar\phi_{h}] = c \cdot I\,
\end{equation}
is satisfied on $X\backslash S$ for some constant $c$, where $I$ is the identity endomorphism of $E$. It is important to note here that, in contrast 
to an admissible metric, the condition defining an admissible Hermitian-Yang-Mills metric depends on the Higgs field. \\

Let $\mathfrak E$ be a torsion-free Higgs sheaf over a compact K\"ahler manifold $X$. Since its singularity set $S$ is a complex analytic subset with 
codimension greater or equal than two, $X\backslash S$ satisfies all assumptions Simpson \cite{Simpson} imposes\footnote{Notice 
that since $X$ is compact, by \cite{Simpson}, Proposition 2.1, it satisfies the assumptions on the base manifold that Simpson introduced. Now, from 
this and \cite{Simpson}, Proposition 2.2, it follows that $X\backslash S$ also satisfies the assumptions.} on the base manifold and hence we can see 
${\mathfrak E}|_{X\backslash S}$ as a Higgs bundle over the non-compact K\"ahler manifold $X\backslash S$. \\

Following Simpson \cite{Simpson}, Proposition 3.3 (see also \cite{Biswas-Schumacher}, Corollary 3.5) it follows that a torsion-free Higgs sheaf over a compact K\"ahler 
manifold with an admissible Hermitian-Yang-Mills metric must be at least semistable. However, as Biswas and Schumacher have shown in \cite{Biswas-Schumacher}, 
this is just a part of a stronger result. To be more precise, they proved the Hitchin-Kobayashi correspondence for Higgs sheaves. This result can be 
written as
\begin{thm}
 Let $\mathfrak E$ be a torsion-free Higgs sheaf over a compact K\"ahler manifold $X$ with K\"ahler form $\omega$. Then, it is $\omega$-polystable if 
and only if there exists an admissible Hermitian-Yang-Mills structure on it.  
\end{thm}

Let $h$ be an admissible metric on $\mathfrak E$, then $K_{h}$ is $L^{1}$-bounded. On the other hand, from \cite{Biswas-Schumacher}, Lemma 2.6, 
we know the Higgs field $\phi$ is also $L^{1}$-bounded on $X\backslash S$ (in particular it is square integrable). From this we conclude that  
\begin{equation}
{\cal K}_{h}=K_{h} + i\Lambda[\phi,\bar\phi_{h}]
\end{equation}
is $L^{1}$-bounded and hence, for any admissible metric $h$ on the torsion-free Higgs sheaf $\mathfrak E$ we must have  
\begin{equation}
 \int_{X\backslash S}|{\cal K}_{h}|\,\omega^{n} < \infty \,. \label{Main L1 property}
\end{equation}

Let ${\rm Herm}^{+}({\mathfrak E}_{X\backslash S})$ be the space of all smooth metrics on ${\mathfrak E}_{X\backslash S}$ satisfying the condition
(\ref{Main L1 property}) and suppose that $h$ and $k$ are two metrics in the same connected component of 
${\rm Herm}^{+}({\mathfrak E}_{X\backslash S})$. Then $h=ke^{v}$ for some endomorphism $v$ of $E|_{X\backslash  S}$ and following Simpson 
\cite{Simpson}, we can write the Donaldson functional as


\begin{equation}
{\cal L}(ke^{v},k)=\int_{X\backslash S}{\rm tr}\left[v({\cal K}_{k}-cI)\right]\frac{\omega^{n}}{n!} + \int_{X\backslash S}\sum_{i,j=1}^{r}\psi(\lambda_{i},\lambda_{j})|{\cal D}''v_{j}^{i}|^{2}\frac{\omega^{n-1}}{(n-1)!} \label{Donaldson functional sheaves} 
\end{equation}
where the function $\psi$ is given by 
\begin{equation}
 \psi(\lambda_{i},\lambda_{j})=\frac{e^{(\lambda_{i}-\lambda_{j})}-(\lambda_{i}-\lambda_{j})-1}{(\lambda_{i}-\lambda_{j})^{2}}\,.
\end{equation}
We define the {\it Donaldson functional on the Higgs sheaf} $\mathfrak E$ just as the corresponding functional (\ref{Donaldson functional sheaves}) 
defined on the Higgs bundle ${\mathfrak E}|_{X\backslash S}$. In \cite{Simpson}, Simpson established an inequality between the supremum of the norm 
of the endomorphism $v$ relating the metrics $h$ and $k$ and the Donaldson functional for Higgs bundles over (non necessarily) compact K\"ahler manifolds; this result can be immediately 
adapted to Higgs sheaves as follows:
\begin{cor}
Let $k$ be an admissible metric on a torsion-free Higgs sheaf $\mathfrak E$ over a compact K\"ahler manifold $X$ with K\"ahler form $\omega$ and 
suppose that ${\rm sup}_{X\backslash S}|{\mathcal K}_{k}|\le B$ for certain fixed constant $B$. If $\mathfrak E$ is $\omega$-stable, then there exist constants 
$C_{1}$ and $C_{2}$ such that
\begin{equation}
 {\rm sup}_{X\backslash S}|v| \le C_{1} + C_{2}\,{\cal L}(ke^{v},k)
\end{equation}
for any selfadjoint endomorphism $v$ with ${\rm tr\,}v=0$ and ${\rm sup}_{X\backslash S}|v|<\infty$ and such that ${\rm sup}_{X\backslash S}|{\cal K}_{ke^{v}}|\le B$.  
\end{cor}

In a future work, we will study more in detail admissible metrics and Donaldson's functional for torsion-free Higgs sheaves and 
the equivalence between semistability and the existence of approximate Hermitian-Yang-Mills metrics for Higgs bundles in higher 
dimensions.

\end{document}